\magnification \magstep1
\input epsf.tex
%\baselineskip .6cm
%\vsize=12truecm
%\hsize=21truecm
\overfullrule=0pt
% Gothic fonts from AMSTeX 
\font\tengoth=eufm10  \font\fivegoth=eufm5
\font\sevengoth=eufm7
\newfam\gothfam  \scriptscriptfont\gothfam=\fivegoth 
\textfont\gothfam=\tengoth \scriptfont\gothfam=\sevengoth

%
% Bold italic fonts 
\font\tenbi=cmmib10  \font\fivebi=cmmib5
\font\sevenbi=cmmib7
\newfam\bifam  \scriptscriptfont\bifam=\fivebi 
\textfont\bifam=\tenbi \scriptfont\bifam=\sevenbi

\font\hd=cmbx10 scaled\magstep1
\def \Box {\hfill\hbox{}\nobreak \vrule width 1.6mm height 1.6mm
depth 0mm  \par \goodbreak \smallskip}

\def \deg  {\mathop{\rm deg}}

\def \Ker {\mathop{\rm Ker}}

\def \iso {\cong}

\def \Hom {\mathop{\rm Hom}}
\def \Tor {{\rm Tor}}
\def \Ext {{\rm Ext}}
\def \H {\mathop {\rm H}\nolimits}
\def \ZZ {{\bf Z}}

\def\dim{\rm dim}

%%
% A small macro producing a box: #1 is the size of the
% box,  #2 is what's in it !
%%

\bigskip
\bigskip

\centerline{\hd Local Cohomology at Monomial Ideals}

\medskip

\centerline{by}

\smallskip

\centerline{\bf Mircea Musta\c{t}\v{a}}

\bigskip
\bigskip
\bigskip

\beginsection {Introduction}

\bigskip

Let $B$ be an ideal in a polynomial ring $R=k[X_1,\ldots, X_n]$
in $n$ variables over a field $k$. The local cohomology of $R$ at $B$
is defined by 
$$\H^i_B(R)={\rm lim}\,\Ext^i_R(R/B^d, R).$$
In general, this limit is not well behaved: the natural maps
$$\Ext^i_R(R/B^d, R)\longrightarrow\H^i_B(R)$$
are not injective and it is difficult to understand how their images
converge to $\H^i_B(R)$ (see Eisenbud, Musta\c{t}\v{a} and Stillman
[1998] for a discussion of related problems).

However, in the case when $B$ is a monomial ideal we will see that
the situation is especially nice if instead of the sequence
$\{B^d\}_{d\geq 1}$ we consider the cofinal sequence of ideals
$\{B_0^{[d]}\}_{d\geq 1}$, consisting of the ``Frobenius powers''
 of the ideal $B_0={\rm radical}(B)$. They are defined as follows:
if $m_1,\ldots ,m_r$ are monomial generators of $B_0$,
then $$B_0^{[d]}=(m_1^d,\ldots,m_r^d).$$

Our first main result is that the natural map
$$\Ext^i_R(R/B_0^{[d]}, R)\longrightarrow\H^i_B(R)$$
is an isomorphism onto the submodule of $H^i_B(R)$
of elements of multidegree $\alpha$, with $\alpha_j\geq -d$ for all $j$. 

The second main result gives a filtration of $\Ext^i_R(R/B, R)$
for a squarefree monomial ideal $B$. For $\alpha\in\{0,1\}^n$,
let ${\rm supp}(\alpha)=\{j\,\vert\,\alpha_j=1\}$ and
$P_{\alpha}=(X_j\,\vert\,j\in {\rm supp}(\alpha))$.

We describe a canonical filtration of
$$\Ext^i_R(R/B, R): 0=M_0\subset\ldots\subset M_n=\Ext^i_R(R/B, R)$$
such that for every $l$,
$$M_l/M_{l-1}\iso\bigoplus_{|\alpha|=l}(R/P_{\alpha}(\alpha))^
{\beta_{l-i,\alpha}(B^{\vee})}.$$
The numbers $\beta_{l-i,\alpha}(B^{\vee})$ are the Betti numbers of
$B^{\vee}$, the Alexander dual ideal of $B$
(see section 3 below for the related definitions). For an interpretation
of this filtration in terms of Betti diagrams, see Remark 1 after
Theorem 3.3, below.

In a slightly weaker form, this result has been conjectured
by David Eisenbud.

\smallskip

Let's see this filtration for a simple example:
$R=k[a,b,c,d]$, $B=(ab, cd)$ and $i=2$.
Since $B$ is a complete intersection,
we get $\Ext^2_R(R/B, R)\iso R/B(1,1,1,1)$.
Our filtration is $M_0=M_1=0$, 
$M_2=R\overline a\overline c+R\overline a\overline d
+R\overline b\overline c+R\overline b\overline d$,
$M_3=R\overline a+R\overline b+R\overline c+R\overline d$
and $M_4=\Ext^2_R(R/B, R)$.

From the description of $\Ext^2_R(R/B, R)$ it follows that
$$M_2/M_1=M_2\iso R/(b, d)(0,1,0,1)\oplus R/(b, c)(0,1,1,0)
\oplus$$
$$ R/(a,d)(1,0,0,1)\oplus R/(a,c)(1,0,1,0),$$
\smallskip
$$M_3/M_2\iso R/(b,c,d)(0,1,1,1)\oplus R/(a,c,d)(1,0,1,1)
\oplus$$
$$ R/(a,b,d)(1,1,0,1)\oplus R/(a,b,c)(1,1,1,0),$$
\smallskip
$$M_4/M_3\iso R/(a,b,c,d)(1,1,1,1).$$

On the other hand,
$B^{\vee}=(bd, bc, ad, ac)$. If $F_{\bullet}$ is the
minimal multigraded resolution of $B^{\vee}$, then
$$F_0=R(0,-1,0,-1)\oplus R(0,-1,-1,0)\oplus R(-1,0,0,-1)
\oplus R(-1,0,-1,0),$$
$$F_1=R(0,-1,-1,-1)\oplus R(-1,0,-1,-1)\oplus R(-1,-1,0,-1)
\oplus R(-1,-1,-1,0),$$
$$F_2=R(-1,-1,-1,-1).$$

We see that for each $\alpha\in\{0,1\}^4$ such that 
$R(-\alpha)$ appears in $F_{l-2}$, there is a corresponding
summand $R/P_{\alpha}(\alpha)$ in $M_l/M_{l-1}$.

\smallskip

In order to prove this result about the filtration of 
$\Ext^i_R(R/B, R)$ we will study the multigraded components
of this module and how an element of the form $X_j\in R$
acts on these components.
As we have seen, it is enough to study
the same problem for $H^i_B(R)$.

We give two descriptions for the degree $\alpha$ part
of $H^i_B(R)$, as simplicial cohomology groups of
certain simplicial complexes depending only on $B$ and
the signs of the components of $\alpha$. The first complex
is on the set of minimal generators of $B$ and the second one
is a full subcomplex of the simplicial complex associated
to $B^{\vee}$ via the Stanley-Reisner correspondence.
The module structure on $H^i_B(R)$ is described by the maps
induced in cohomology by inclusion of simplicial complexes. 

As a first consequence of these results and using also
a formula of Hochster [1977], we obtain an isomorphism 
$$\Ext^i_R(R/B, R)_{-\alpha}\iso\Tor^R_{|\alpha|-i}(B^{\vee}, k)_{\alpha},$$
for every $\alpha\in\{0,1\}^n$.

This result is equivalent to the fact that in our filtration the numbers
are as stated above. This isomorphism has been obtained also
by Yanagawa [1998]. It can be considerd as a strong form of the
inequality of Bayer, Charalambous and Popescu [1998] between
the Betti numbers of $B$ and those of $B^{\vee}$. As shown in that
paper, this implies that $B$ and $B^{\vee}$ have the same
extremal Betti numbers, extending results of Eagon and Reiner [1996]
and Terai [1997].  

As a final application of our analysis of the graded pieces of $\Ext^i_R
(R/B, R)$, we give a topological description for
the associated primes of $\Ext^i_R(R/B, R)$.
In the terminology of Vasconcelos [1998], these are the homological
associated primes of $R/B$.
 In particular, we characterize
the minimal associated primes of $\Ext^i_R(R/B, R)$
 using only the Betti numbers of $B^{\vee}$.

\smallskip

We mention here the recent work of Terai [1998]
on the Hilbert function of the modules $H^i_B(R)$. 
It is easy to see that using the results in our paper
one can deduce Terai's formula for this Hilbert function.

The problem of effectively computing the local cohomology modules
with respect to an arbitrary ideal is quite
difficult since these modules are not finitely generated.
The general approach is to use the $D$-module structure for
the local cohomology (see, for example, Walther [1999]).
However, in the special case of monomial ideals our results show
that it is possible to make this computation with elementary methods. 

Our main motivation for studying local cohomology at monomial ideals
comes from the applications in the context of toric varieties. Via
the homogeneous coordinate ring, the cohomology of sheaves on such
a variety can be expressed as local cohomology of modules at the
``irrelevant ideal'', which is a squarefree monomial ideal. For a method
of computing the cohomology of sheaves on toric varieties 
in this way, see Eisenbud,
Musta\c{t}\v{a} and Stillman [1998]. For applications to vanishing theorems on 
toric varieties and related results, see Musta\c{t}\v{a} [1999].

\smallskip

The main reference for the definitions and the results that we use is 
Eisenbud [1995].
For the basic facts about the cohomology of simplicial complexes, see 
Munkres [1984].
Cohomology of simplicial complexes is always taken to be reduced
cohomology. Notice also that we make a distinction
between the empty complex which contains just the empty set (which has
nontrivial cohomology in degree $-1$) and the void complex which 
doesn't contain any set (whose cohomology is trivial in any degree).

\smallskip

This work has been done in connection with a joint project with
David Eisenbud and Mike Stillman. We would like to thank them for their
constant encouragement and for generously sharing their insight with us.
We are also grateful to Josep Alvarez Montaner for pointing out 
a mistake in an earlier version of this paper.

\bigskip\bigskip

\beginsection \S 1. Local cohomology as a union of Ext modules
  
\bigskip

Let $B\subset R=k[X_1,\ldots ,X_n]$ be a squarefree monomial ideal.
 All the modules which 
appear are $\ZZ^n$-graded.
We partially order the elements of $\ZZ^n$
by setting $\alpha\geq\beta$ iff $\alpha_j\geq
\beta_j$, for all $j$. 

\smallskip

\proclaim Theorem 1.1. For each $i$ and $d$, the natural map
$$\Ext^i_R(R/B^{[d]},R)\longrightarrow\H^i_B(R)$$
is an isomorphism onto the submodule of $\H^i_B(R)$
of elements of degree $\geq (-d,\ldots, -d)$.

\smallskip

\noindent{\sl Proof.\ } We will compute $\Ext^i_R(R/B^{[d]}, R)$
using the Taylor resolution $F_{\bullet}^d$ of $R/B^{[d]}$ ( see 
Eisenbud [1995], exercise 17.11).
 The inclusion $B^{[d+1]}\longrightarrow B^{[d]}$,
$d\geq 1$ induces a morphism of complexes $\phi^d : F^{d+1}_{\bullet}
\longrightarrow F^d_{\bullet}$. The assertions in the theorem 
are consequences of the more precise lemma below.

\smallskip

 \proclaim Lemma 1.2. If $(\phi^d)^* : (F^d_{\bullet})^*
\longrightarrow(F^{d+1}_{\bullet})^*$ is the dual
$\Hom_R(\phi^d, R)$ of the above map, then in a multidegree $\alpha\in\ZZ^n$
we have:
\item {(a)} If $\alpha\geq (-d,\ldots, -d)$, then $(\phi^d)^*_{\alpha}$
is an isomorphism of complexes.
\item {(b)} If $\alpha_j<-d$ for some $j$, $1\leq j\leq n$, then
$(F^d_{\bullet})^*_{\alpha}=0$, so $(\phi^d)^*_{\alpha}$
is the zero map.

\smallskip

\noindent {\sl Proof of the lemma.\ }                            
Let $m_1,\ldots, m_r$ be monomial minimal generators of $B$.
 For any subset $I$ of
$\{1,\dots,r\}$ we set 
$$m_I = {\rm LCM}\,\{m_i\mid i\in I\}.$$
As each $m_I$ is square-free, $\deg m_I\in \ZZ^n$ is a vector
of ones and zeros.

Recall from Eisenbud [1995] that $F^d_{\bullet}$
 is a free $R$-module with basis
$\{f^d_I\,|\,I\subset\{1,\ldots, r\}\}$, where $\deg(f^d_I)=d\,\deg(m_I)$.
Therefore, the degree $\alpha$ part of
 $(F^d_{\bullet})^*$ has a vector space basis
consisting of elements of the form
$n e^d_I$ where $n\in R$ is a monomial, $e^d_I=(f^d_I)^*$ 
has degree equal to $-d\deg(m_I)$, and 
$
\deg(n)-d\deg(m_I)=\alpha.
$

Part (b) of the Lemma follows at once. For part (a), note that
$(\phi^k)^* : (F^d_{\bullet})^* \longrightarrow    
 (F^{d+1}_{\bullet})^*$ takes $e^d_I$ to
$m_Ie^{d+1}_I$. The vector $\deg(e^{d+1}_I)=-(d+1)\deg(m_I)$ has
entry $-(d+1)$ wherever $\deg(m_I)$ has entry 1, so
any element $n e^{d+1}_I$ of degree 
$\alpha\geq (-d,-d,\dots,-d)$ 
must have $n$ divisible by $m_I$. It is thus of the form
$(\phi^d)^*(x)$ for the unique element $x=(n/m_I)e^d_I$,
as required.\Box

\bigskip

\beginsection \S 2. Local cohomology as simplicial cohomology

\bigskip

To describe $\H^i_B(R)$ in a multidegree $\alpha\in\ZZ^n$, we will
use two simplicial complexes associated with $B$ and
$\alpha$. We will assume that $B\neq (0)$.

By computing local cohomology using the Taylor complex
we will express $H^i_B(R)_{\alpha}$ as the simplicial cohomology
of a complex
on the set of minimal generators of $B$. We will interpret
this later as the cohomology of an other complex, this time
on the potentially smaller set $\{1,\ldots, n\}$.  
This one is a full subcomplex of the complex associated to the
dual ideal $B^{\vee}$ via the Stanley-Reisner correspondence.
In fact, this is the complex used in the computation of the Betti numbers
of $B^{\vee}$ (see the next section for the definitions).
We will use this result to derive the relation between
$\Ext(R/B, R)$ and $\Tor^R(B^{\vee}, k)$ in Corolary 3.1 below. 

\smallskip

Let $m_1,\ldots, m_r$ be the minimal monomial generators of $B$.
As above, for $J\subset\{1,\ldots, r\}$, $m_J$ will denote
${\rm LCM}\,(m_j; j\in J)$.

For $i\in\{1,\ldots, n\}$, we define
$$T_i:=\{J\subset\{1,\ldots, r\}\,|\, X_i\not |\, m_J\}.$$
For every subset $I\subset\{1,\ldots, n\}$, we define 
$T_I=\bigcup_{i\in I} T_i$.
When $I=\emptyset$, we take $T_I$ to be the void complex. It is clear
that each $T_i$ is a simplicial complex on the set $\{1,\ldots, r\}$,
and therefore so is $T_I$.

For $\alpha\in\ZZ^n$, we take $I_{\alpha}=\{i\,|\,\alpha_i\leq-1\}\subset
\{1,\ldots, n\}$. Note that the complex $T_{I_{\alpha}}$ depends only on
the signs of the components of $\alpha$ (and, of course, on $B$).

If $e_1,\ldots,e_n$ is the canonical basis of $\ZZ^n$ and
 $\alpha'=\alpha+e_l$, we have obviously $I_{\alpha'}\subset
I_{\alpha}$, with equality iff $\alpha_l\neq -1$. Therefore,
$T_{I_{\alpha'}}$ is a subcomplex of $T_{I_{\alpha}}$.

\smallskip

\proclaim Theorem 2.1.
\item {(a)} With the above notation, we have 
$$\H_B^i(R)_{\alpha}\iso \H^{i-2}(T_{I_{\alpha}}; k).$$

\item {(b)} Via the isomorphisms given in (a), the multiplication 
by $X_l$:
$$\nu_{X_l}\,:\,\H^i_B(R)_{\alpha}\longrightarrow
\H^i_B(R)_{\alpha'}$$
corresponds to the morphism:
$$\H^{i-2}(T_{I_{\alpha}}; k)\longrightarrow\H^{i-2}(T_{I_{\alpha'}}; k),$$
induced in cohomology by the inclusion $T_{I_{\alpha'}}\subset
T_{I_{\alpha}}$. In particular, if $\alpha_l\neq -1$, then
$\nu_{X_l}$ is an isomorphism.

\smallskip

\noindent{\sl Proof.\ } We have seen in Lemma 1.2 that
$$\Ext^i_R(R/B^{[d]}, R)_{\alpha}\iso\H_B^i(R)_{\alpha}$$
if $\alpha\geq (-d,\ldots, -d)$. We fix such a $d$.
With the notations in Lemma 1.2~, we have seen that
the degree $\alpha$ part of 
$(F^d_{\bullet})^*$ has a vector space
 basis consisting of elements of the form 
$n e_J^d$, where $n\in R$ is a monomial and $\deg(n)-d\deg(m_J)=\alpha$.
Therefore, the basis
of $(F^d_p)^*_{\alpha}$ is indexed by those $J\subset\{1,\ldots ,r\}$ with
$|J|=p$ and $\alpha+d\deg(m_J)\geq (0,\ldots ,0)$.
Because $\alpha_j\leq -1$ iff $j\in I_{\alpha}$
 and $\alpha\geq (-d,\ldots ,-d)$, the above inequality
is equivalent to $X_j\vert m_J$ for every $j\in I_{\alpha}$
i.e. to $J\not\in T_{I_{\alpha}}$.

Let $G^{\bullet}$ be the cochain complex computing the relative cohomology
of the pair $(D, T_{I_{\alpha}})$ with 
coefficients in $k$, where $D$ is the full simplicial complex
on the set $\{1,\ldots, r\}$.

If $I_{\alpha}\neq\emptyset$,
 then the degree $\alpha$ part of $(F^d_p)^*$ is equal
to $G^{p-1}$ for every $p$. Moreover, the maps are the same and therefore
we get $H_B^i(R)_{\alpha}\iso  
\H^{i-1}(D,T_{I_{\alpha}}; k)$. Since $D$ is contractible, 
the long exact sequence in cohomology of the pair $(D, T_{I_{\alpha}})$ yields
$\H_B^i(R)_{\alpha}\iso \H^{i-2}(T_{I_{\alpha}}; k)$.

If $I_{\alpha}=\emptyset$,
 then $(F^d_{\bullet})^*$ in degree $\alpha$ is up to a shift the complex
computing the reduced cohomology of $D$ with 
coefficients in $k$. Since $D$ is contractible, we get
 $\H_B^i(R)_{\alpha}=0=\H^{i-2}(T_{I_{\alpha}}; k)$,
which completes the proof of part (a).

For part (b), we may suppose that $I_{\alpha'}\neq (0)$.
With the above notations, $\nu_{X_l}$ is induced by the map
$\phi_l\,:\,(F^d_p)^*_{\alpha}\longrightarrow
(F^d_p)^*_{\alpha'}$, given by $\phi_l(n e^d_J)=X_l n e^d_J$.

If $G'^{\bullet}$ is constructed as above, but for $\alpha'$ instead
of $\alpha$, then via the isomorphisms:
$$(F^d_p)^*_{\alpha}\iso G^{p-1},$$
$$(F^d_p)^*_{\alpha'}\iso G'^{p-1},$$
the map $\phi_l$ corresponds to the canonical projection
$G^{p-1}\longrightarrow G'^{p-1}$, which concludes the proof of part (b).\Box

\medskip

\proclaim Remark.
The last assertion in Theorem 2.1(b), that $\nu_{X_l}$
is an isomorphism if $\alpha_l\neq -1$ has been obtained
also in Yanagawa [1998].

\medskip

The next corollary describes $\H^i_B(R)_{\alpha}$ as the cohomology of a 
simplicial complex with vertex set $\{1,\ldots, n\}$. 

 We first introduce 
the complex $\Delta$ defined by:
$$ \Delta:= \{F\subset\{1,\ldots, n\}\,|\,\prod_{j\not\in F}X_j\in B\}.$$
In fact , by the Stanley-Reisner correspondence between square-free
monomial ideals and simplicial complexes (see 
Bruns and Herzog [1993]), $\Delta$ corresponds
to $B^{\vee}$.

For any subset $I\subset\{1,\ldots, n\}$, we define $\Delta_I$ to be the
full simplicial subcomplex of $\Delta$ supported on $I$:
$$\Delta_I:=\{F\subset\{1,\ldots, n\}\,|\,F\in\Delta , F\subset I\}.$$ 
When $I=\emptyset$, we take $\Delta_I$ to be the void complex. It is clear
that if $I\subset I'$, then $\Delta_{I'}$ is a subcomplex
of $\Delta_I$. This is the case if $\alpha'=\alpha+e_l$,
$I=I_{\alpha}$ and $I'=I_{\alpha'}$.

\smallskip 

\proclaim Corollary 2.2.
\item {(a)} With the above notation, for any $\alpha\in\ZZ^n$
$$\H^i_B(R)_{\alpha}\iso\H^{i-2}(\Delta_{I_{\alpha}}; k).$$ 

\item {(b)} Via the isomorphisms given by (a), the multiplication map
$\nu_{X_l}$ corresponds to the morphism:
$$\H^{i-2}(\Delta_{I_{\alpha}}; k)\longrightarrow
\H^{i-2}(\Delta_{I_{\alpha'}}; k),$$
induced in cohomology by the inclusion $\Delta_{I_{\alpha'}}
\subset\Delta_{I_{\alpha}}$.
\smallskip

\noindent {\sl Proof.\ } Using the notation in Theorem 2.1,  
 if $I_{\alpha}\neq\emptyset$,
 then $T_{I_{\alpha}}=\bigcup_{i\in I_{\alpha}}T_i$.

If $i_1,\ldots, i_k\in I_{\alpha}$ and $\bigcap_{1\leq p\leq k}T_{i_p}
\neq\emptyset$, then
$$\bigcap_{1\leq p\leq k}T_{i_p}=
\{J\subset\{1,\ldots, r\}\,|\,X_{i_p}\not |\,m_J, 1\leq p\leq k\}$$  
is the full simplicial complex on those $j$ with $X_{i_p}\not |\,m_j$,
for every $p$, $1\leq p\leq k$. Therefore it is contractible.

This shows that we can compute the cohomology of $T_I$ as the cohomology
of the nerve ${\cal N}$ of the cover $T_I=\bigcup_{i\in I}T_i$ (see 
Godement [1958]).
But by definition, $\{i_1,\ldots, i_k\}\subset I$ is a simplex in ${\cal N}$
iff $\bigcap_{1\leq p\leq k}T_{i_p}\neq\emptyset$ iff
there is $j$ such that $X_{i_p}\not |\,m_j$ for every p,
$1\leq p\leq k$. This shows that ${\cal N}=\Delta_I$ and we get that
$\H_B^i(R)_{\alpha}\iso\H^{i-2}(\Delta_I; k)$ when $I\neq\emptyset$.

When $I=\emptyset$, $\H^i_B(R)_{\alpha}=0$ by theorem 2.1 and also 
$\H^{i-2}(\Delta_I; k)=0$ (the reduced cohomology
of the void simplicial complex is zero).

Part (b) follows immediately from part (b) in Theorem 2.1
and the fact that the isomorphism between the cohomology of a space
and that of the nerve of a cover as above is functorial. \Box

\bigskip

\proclaim Remark. The same type of arguments as in the proofs of 
Theorem 2.1 and of Corollary 2.2 can be used to give a topological
description for $\Ext^i_B(R/B, R)_{\alpha}$, for a possibly
non-reduced nonzero monomial ideal $B$. 
Namely, for $\alpha\in\ZZ^n$, we define the simplicial complex 
$\Delta_{\alpha}$
on $\{1,\ldots, n\}$ by $J\in\Delta_{\alpha}$ iff there is a
monomial $m$ in $B$ such that $\deg(X^{\alpha} m)_j<0$  
for $j\in J$. We make the convention that $\Delta$ is the void
complex iff $\alpha\geq 0$. Then $$ \Ext^i_R(R/B, R)_{\alpha}\iso
\H^{i-2}(\Delta_{\alpha}; k).$$

Moreover, we can describe these $k$-vector spaces using a more
geometric object. If we view $B\subset\ZZ^n\subset\,{\bf R}^n$, let
$P_{\alpha}$ be the subspace of ${\rm R}^n$ supported on $B$,
translated by $\alpha$, minus the first quadrant. More precisely,
$$P_{\alpha}=\{x\in\,{\bf R}^n | x-\alpha\geq m,\,{\rm for}\,
{\rm some}\,m\in B\}\setminus{\bf R}_+^n.$$
Then, using a similar argument to the one in the proof of corollary 1.4,
one can show that $$\Ext^i_R(R/B, R)_{\alpha}\iso\H^{i-2}(P_{\alpha}; k),$$
where the right-hand side is the reduced singular cohomology group.
Here we have to make the convention that for $\alpha\geq 0$,
$P_{\alpha}$ is the ``void topological space'', with trivial reduced
cohomology (as oposed to the empty topological space which has
nonzero reduced cohomology in degree $-1$).

We leave the details of the proof to the interested reader. 

\bigskip

\beginsection \S 3. The filtration on the Ext modules

\bigskip

The Alexander dual of a reduced monomial ideal $B$ is defined by
$$B^{\vee}=(X^F | F\subset\{1,\ldots n\}, X^{F^c}\notin B),$$
 where $F^c:=\{1,\ldots, n\}\setminus F$ (see Bayer, Charalambous
and Popescu [1998] for interpretation in terms
of Alexander duality~). Note that $(B^{\vee})^{\vee}=B$.

We will derive first a relation between $\Ext_R(R/B, R)$ and
$\Tor^R(B^{\vee}, k)$. This can be seen as a stronger form
of the inequality in Bayer, Charalambous and Popescu [1998]
 between the Betti numbers of $B$ and $B^{\vee}$.
  
For $\alpha\in\ZZ^n$, we will denote $|\alpha|=\sum_i\alpha_i$.

\smallskip

 \proclaim Corollary 3.1. Let $B\subset R=k[X_1,\ldots, X_n]$ be a reduced 
monomial ideal and $\alpha\in\ZZ^n$ a multidegree. If $\alpha\notin\{0,1\}^n$,
then $\Tor^R_i(B^{\vee}, k)_{\alpha}=0$, and if 
$\alpha\in\{0,1\}^n$, then
$$\Tor^R_i(B^{\vee}, k)_{\alpha}\iso\Ext_R^{|\alpha|-i}(R/B, R)_
{-\alpha}.$$

\smallskip

\noindent{\sl Proof.\ } We will use Hochster's formula for
 the Betti numbers of reduced
monomial ideals (see, for example, 
Hochster [1977] or Bayer, Charalambous and Popescu [1998]). It says that 
if ${\alpha}\notin\{0,1\}^n$, then
$\Tor^R_i(B^{\vee}, k)_{\alpha}=0$  and if 
${\alpha}\in\{0,1\}^n$, then 
$$\Tor^R_i(B^{\vee}, k)_{\alpha}\iso \H^{|\alpha|-i-2}(\Delta_I; k),$$
where $I$ is the support of $\alpha$.

Obviously, we may suppose that $B\neq (0)$. If $\alpha\in\{0,1\}^n$, then 
corollary 2.2 gives
 $$\H^{|\alpha|-i-2}(\Delta_I; k)\iso\H_B^{|\alpha|-i}(R)_{-\alpha}$$
and theorem 1.1 gives $$\H_B^{|\alpha|-i}(R)_{-\alpha}
\iso\Ext_R^{|\alpha|-i}(R/B, R)_{-\alpha}.$$
Putting together these isomorphisms, we get the assertion of the 
corollary. \Box

\medskip

We recall that the multigraded Betti numbers of $B$ are defined by 
$$\beta_{i,\alpha}(B):=\dim_k \Tor_i^R(B, k)_{\alpha}.$$ 
Equivalently, if $F_{\bullet}$ is a multigraded 
minimal resolution of $B$, then $$F_i\iso\sum_{\alpha\in\ZZ^n}
R(-\alpha)^{\beta_{i,\alpha}(B)}.$$

One says that $(i,\alpha)$ is extremal (or that $\beta_{i,\alpha}$ is
extremal) if $\beta_{j,\alpha'}(B)=0$ for all $j\geq i$ and
 $\alpha'>\alpha$
such that $|\alpha'|-|\alpha|\geq j-i$. 

\medskip

\proclaim Remark. Using Theorems 1.1, 2.1(b) and Corollary 3.1
one can give a formula for the Hilbert function
of $\H^i_B(R)$
using the Betti numbers of $B^{\vee}$. This formula is 
equivalent to the one which appears in Terai [1998].

\medskip

As a consequence of the above corollary, we obtain the inequality 
between the Betti numbers of $B$ and $B^{\vee}$ from 
Bayer, Charalambous and Popescu [1998].
It implies the equality of extremal Betti numbers 
from that paper, in particular
the equality ${\rm reg}\,B={\rm pd}(R/B^{\vee})$ 
from Terai [1997].

\proclaim Corollary 3.2. If $B\subset R$ is a reduced monomial ideal, then
$$\beta_{i,\alpha}(B)\leq
\sum_{\alpha\leq\alpha' }\beta_{|\alpha|-i-1,\alpha'}(B^{\vee}),$$
for every $i\geq 0$ and every $\alpha\in\{0,1\}^n$.
If $\beta_{|\alpha|-i-1,\alpha}(B^{\vee})$ is extremal, then so is 
$\beta_{i,\alpha}(B)$ and 
$$\beta_{i,\alpha}(B)=\beta_{|\alpha|-i-1,\alpha}(B^{\vee}).$$

\noindent {\sl Proof.\ } Since $\beta_{i,\alpha}(B)
=\dim_k \Tor_i^R(B, k)_{\alpha}$, by the previous
corollary we get
$$\beta_{i,\alpha}(B)=
\dim_k\Ext_R^{|\alpha|-i}(R/B^{\vee}, R)_{-\alpha}=
\dim_k\H^{|\alpha|-i}({\rm Hom}(F_{\bullet}, R))_{-\alpha},$$
where $F_{\bullet}$ is the minimal free resolution of $R/B^{\vee}$.

Since $F_{|\alpha|-i}=
\oplus_{\alpha'\in\ZZ^n}R(-\alpha')^{\beta_{|\alpha|-i-1,\alpha'}(B^{\vee})}$,
we get $$\beta_{i,\alpha}(B)
\leq\sum_{\alpha'\in\ZZ^n}\beta_{|\alpha|-i-1,
\alpha'}(B^{\vee})\dim_k(R(\alpha')_{-\alpha})=
\sum_{\alpha\leq\alpha'}\beta_{|\alpha|-i-1,\alpha'}(B^{\vee}).$$

If $\beta_{|\alpha|-i-1,\alpha}
(B^{\vee})$ is extremal, the above inequality
becomes $\beta_{i,\alpha}(B)\leq\beta_{|\alpha|-i-1,\alpha}(B^{\vee})$.
Applying the same inequality for $j\geq i$ and $\alpha'>\alpha$ such that 
$|\alpha'|-|\alpha|\geq j-i$ and the fact that
 $\beta_{|\alpha|-i-1,\alpha}(B^{\vee})$ is extremal,
we get that $\beta_{i,\alpha}(B)$ is extremal.

Applying the previous inequality with $B$ replaced by $B^{\vee}$, we obtain
$\beta_{|\alpha|-i-1,\alpha}(B^{\vee})
\leq\beta_{i,\alpha}(B)$, which concludes the proof.\Box

\bigskip
We fix some notations for the remaining of this section.
Let $[n]=\{0,1\}^n$ and $[n]_l=\{\alpha\in [n]\,\vert\, |\alpha|=l\}$,
for every $l$, $0\leq l\leq n$.
For $\alpha\in [n]$, let ${\rm supp}(\alpha)=\{j\,\vert\,
\alpha_j=1\}$ and
 $P_{\alpha}=(X_j\,\vert\,j\in {\rm supp}(\alpha))$.
 The ideals $P_{\alpha}$, $\alpha\in [n]$
are exactly the monomial prime ideals of $R$. 

The following theorem gives the canonical filtration
of $\Ext^i_R(R/B, R)$ announced in the Introduction.

\proclaim Theorem 3.3. Let $B\subset R$ be a squarefree monomial
ideal. For each $l$, $0\leq l\leq n$, let $M_l$ be the submodule
of $\Ext^i_R(R/B, R)$ generated by all $\Ext^i_R(R/B, R)_{-\alpha}$,
for $\alpha\in [n]$, $|\alpha|\leq l$. Then $M_0=0$,
$M_n=\Ext^i_R(R/B, R)$ and for every $l$, $0\leq l\leq n$,
$$M_l/M_{l-1}\iso\bigoplus_{\alpha\in [n]_l}
(R/P_{\alpha}(\alpha))^{\beta_{l-i,
\alpha}(B^{\vee})}.$$

\noindent {\sl Proof.\ } Clearly we may suppose $B\neq 0$.
The fact that $M_0=0$ follows from Corollary 2.2(a).

Let's see first that $M_n=\Ext^i_R(R/B, R)$.
 For this it is enough to prove that all the minimal monomial
generators of $\Ext^i_R(R/B, R)$ are in degrees $-\alpha$,
$\alpha\in [n]$.

Indeed, if $\alpha_j\leq -1$ for some $j$, then the multiplication
by $X_j$ defines an isomorphism $$\Ext^i_R(R/B, R)_{-\alpha-e_j}
\longrightarrow\Ext^i_R(R/B, R)_{-\alpha}$$
by Corollary 2.2(b) and Theorem 1.1. In particular, there are no minimal
generators in degree $-\alpha$.  

On the other hand, by Theorem 1.1, $\Ext^i_R(R/B, R)_{-\alpha}=0$
if $\alpha_j\geq 2$, for some $j$. Therefore we have $M_n=\Ext^i_R(R/B, R)$.

Suppose now that we have homogeneous elements
$f_1,\ldots, f_r$ with $\deg(f_q)\in [n]_{l'}$, $l'\leq l$,
for every $q$, $1\leq q\leq r$.
We suppose  that they are linearly independent
 over $k$ and that their linear span contains
$\Ext^i_R(R/B, R)_{-\alpha}$, for every $\alpha\in [n]_{l'}$,
$l'\leq l-1$.
We will suppose also that $\deg(f_r)=-\alpha$, $|\alpha|=l$.
If $T:=\sum_{1\leq q\leq r-1}R\,f_q$, let $\overline f_r$ be
the image of $f_r$ in $M_l/T$.
\smallskip

{\sl Claim.} With the above notations, ${\rm Ann}_R(\overline f_r)
=P_{\alpha}$. 

Let $F={\rm supp}(\alpha)$.
If $j\in F$, then $\deg(X_jf_r)=-(\alpha-e_j)$,
$\alpha-e_j\in [n]$.
By our assumption, it follows that $X_jf_r\in T$, so that
$P_{\alpha}\subset{\rm Ann}_R(\overline f_r)$.

Conversely, consider now $m=\prod X_j^{m_j}\in
{\rm Ann}\,\overline f_r$ and suppose that
$m\not\in (X_j\,|\,j\in F)$. We can suppose that $m$ has minimal degree.
 Let $j$ be such that $m_j\geq 1$. Then $j
\not\in F$ and therefore $m_j-\alpha_j=m_j\geq 1$.
Since $m\,f_r\in T$, we can write
$$m\,f_r=\sum_{q<r}c_q n_q f_q,$$
where $n_q$ are monomials and $c_q\in k$. Since $\deg(f_q)\leq 0$
for every $q$, in the above equality we may assume that 
$X_j|n_q$ for every $q$ such that $c_q\neq 0$.
But by Corollary 2.2(b) and Theorem 1.1,
 the multiplication by $X_j$ is an isomorphism:
$$\Ext^i_R(R/B, R)_{-\alpha+\deg\,m-e_j}\longrightarrow
\Ext^i_R(R/B, R)_{-\alpha+\deg\,m}.$$
Therefore $m/X_j\in{\rm Ann}\overline f_r$, in contradiction with the
minimality of $m$. We get ${\rm Ann}\overline f_r
=(X_j\,|\,j\in F)$, which completes the proof of the claim.

\smallskip

The first consequence is that for every nonzero $f\in M_l$,
$\deg(f)=-\alpha$, $\alpha\in [n]_l$, if
 $\overline f$ is the image of $f$
in $M_l/M_{l-1}$, then ${\rm Ann}_R(\overline f)=P_{\alpha}$, so
that $R\,\overline f\iso R/P_{\alpha}(\alpha)$.

Let's consider now a homogeneous basis $f_1,\ldots, f_N$
of $\oplus_{\alpha\in [n]_l}\Ext^i_R(R/B, R)_{-\alpha}$.
By Corollary 3.1,
$$\dim_k\Ext^i_R(R/B, R)_{-\alpha}=\beta_{l-i, \alpha}
(B^{\vee}).$$ 
Therefore, to complete the proof of the theorem, it is enough to
show that
 $$M_l/M_{l-1}\iso\bigoplus_{1\leq j\leq N}R\overline f_j.$$
Here $\overline f_j$ denotes the image of $f_j$ in $M_l/M_{l-1}$.

Since $M_l=M_{l-1}+\sum_{1\leq j\leq N}R f_j$, we have only to show
that if 
$\sum_{1\leq j\leq N}n_jf_j\in M_{l-1}$, then $n_jf_j\in M_{l-1}$
for every $j$, $1\leq j\leq N$.

Let $\{g_1,\ldots, g_{N'}\}$ be the union of
 homogeneous bases for $\Ext^i_R(R/B, R)_{-\alpha}$,
for $\alpha\in [n]_{l'}$, $l'\leq l-1$.

Let's fix some $j$, with $1\leq j\leq N$. If $\deg(f_j)=-\alpha$,
by applying the above claim to $f_j$, as part of 
$\{f_p\,\vert\,1\leq p\leq N\}\cup\{g_{p'}\,\vert\,1\leq p'\leq N'\}$,
we get that $n_j\in P_{\alpha}$. But we have already seen that
$P_{\alpha}f_j\subset M_{l-1}$ and therefore the proof is complete. \Box

\smallskip
\proclaim Remark 1. We can interpret the statement of Theorem 3.3
using the multigraded Betti diagram of $B^{\vee}$.
This is the diagram having at the intersection of the 
$i^{\rm th}$ row with the $j^{\rm th}$ column the Betti numbers
$\beta_{j,\alpha}(B^{\vee})$, for $\alpha\in\ZZ^n$,
 $|\alpha|=i+j$.

For each $i$ and $j$ we form a module corresponding to $(i, j)$:
$$E_{i,j}=\bigoplus_{\alpha\in [n]_{i+j}}(R/P_{\alpha}(\alpha))
^{\beta_{j,\alpha}(B^{\vee})}.$$
Theorem 3.3 gives a filtration of $\Ext^i_R(R/B, R)$ having as
quotients the modules constructed above corresponding to the
$i^{\rm th}$ row: $E_{i, j}$, $j\in\ZZ$.

Notice that by definition, $\Tor^R_i(B^{\vee}, k)$ is obtained
by a ``dual'' procedure applied to the $i^{\rm th}$ column
(in this case the extensions being trivial).
Indeed, if for $(j, i)$ we put
$$E'_{j,i}=\bigoplus_{\alpha\in [n]_{i+j}}k(-\alpha)
^{\beta_{i,\alpha}(B^{\vee})},$$
then $\Tor^R_i(B^{\vee}, k)\iso\oplus_{j\in\ZZ}E'_{j, i}$.

\smallskip
\proclaim Remark 2. Using Theorem 3.3 one can compute the Hilbert series
of $\Ext^i_R(R/B, R)$ in terms of the Betti numbers of $B^{\vee}$.
Using local duality, one can derive the fomula, due to Hochster [1997],
for the Hilbert series of the local cohomology modules
$H^{n-i}_{\underline{m}}(R/B)$, where $\underline{m}=(X_1,\ldots, X_n)$
 (see also Bruns and Herzog [1993], Theorem 5.3.8).

\bigskip

We describe now the set of homological associated primes of $R/B$
i.e. the set $$\cup_{i\geq 0}{\rm Ass}(\Ext^i_R(R/B, R))$$ 
(see Vasconcelos [1998]).
Since the module $\Ext^i_R(R/B, R)$ is $\ZZ^n$-graded, its associated
primes are of the form $P_{\alpha}$, for some $\alpha\in [n]$.
In fact, Theorem 3.3 shows that
$${\rm Ass}(\Ext^i_R(R/B, R))\subset\{P_{\alpha}\,\vert\,
\beta_{|\alpha|-i, \alpha}(B^{\vee})\neq 0\}.$$

The next result gives the necessary and sufficient condition
for a prime ideal $P_{\alpha}$ to be in ${\rm Ass}(\Ext^i_R(R/B, R))$.
In particular, we get the characterization of the minimal associated
primes of this module using only the Betti numbers of $B^{\vee}$.

\proclaim Theorem 3.4. Let $B\subset R$ be a nonzero square-free
monomial ideal and $\alpha\in [n]$. Let $F={\rm supp}(\alpha)$.

\item {(a)} The ideal $P_{\alpha}$ belongs to ${\rm Ass}(\Ext^i_R(R/B, R))$ iff
$$\bigcap_{j\in F}\Ker (H^{i-2}(\Delta_F; k)\longrightarrow
\H^{i-2}(\Delta_{F\setminus j}; k))\neq 0.$$

\item {(b)} The ideal $P_{\alpha}$ is a minimal prime in 
${\rm Ass}(\Ext^i_R(R/B, R))$ iff
$$\beta_{|\alpha|-i, \alpha}(B^{\vee})\neq 0$$
and 
$$\beta_{|\alpha'|-i, \alpha'}(B^{\vee})=0,$$
for every $\alpha '\in [n]$, $\alpha '\leq\alpha$, $\alpha ' \neq\alpha$.

\smallskip

\noindent {\sl Proof.\ }  By Corollary 2.2, the condition in (a)
is equivalent to the existence of $u\in\Ext^i_R(R/B, R)_{-\alpha}$,
$u\neq 0$ such that $X_j u=0$ for every $j\in F$.
Since $\alpha_j=0$ for $j\not\in F$, Corollary 2.2(b) and
Theorem 1.1 imply that 
for every monomial $m$, $m\not\in P_{\alpha}$,
the multiplication by $m$ is injective on
$\Ext^i_R(R/B, R)_{-\alpha}$.

Therefore, in the above situation we have
${\rm Ann}_R(u)=P_{\alpha}$, so that 
 $P_{\alpha}$ is an element of ${\rm Ass}(\Ext^i_R(R/B, R))$.

Conversely, suppose that $P_{\alpha}\in {\rm Ass}(Ext^i_R(R/B, R))$.
Since $P_{\alpha}$
 and $\Ext^i_R(R/B, R)$ are $\ZZ^n$-graded, this is 
equivalent to the existence of $u\in\Ext^i_R(R/B, R)_{\alpha'}$, for some
$\alpha'\in\ZZ^n$, such that $P_{\alpha}={\rm Ann}_R(u)$.
 To complete the proof of
part (a), it is enough to show that we can take $\alpha'=-\alpha$.

By Theorem 1.1, $\alpha'\geq (-1,\ldots, -1)$. Since $X_j u=0$ for $j\in F$,
multiplication by $X_j$ on $\Ext^i_R(R/B, R)_{\alpha}$ is not injective
so that by Corollary 2.2(b), we must have $\alpha'_j=-1$ for $j\in F$.

Let's consider some $j\not\in F$. If $\alpha'_j\geq 1$, by Corollary 2.2(b)
there is $u'\in\Ext^i_R(R/B, R)_{\alpha''}$, $\alpha''=\alpha'-\alpha'_j e_j$
such that $X_j^{\alpha'_j}u'=u$ and ${\rm Ann}_R(u')={\rm Ann}_R(u)
=P_{\alpha}$.
Therefore, we may suppose that $\alpha'_j\leq 0$.

If $\alpha'_j=-1$, since $X_j\not\in{\rm Ann}_R(u)$, which is prime, we have
${\rm Ann}_R(X_j u)={\rm Ann}_R(u)=P_{\alpha}$. This shows that we may suppose
$\alpha'_j=0$ for every $j\not\in F$, so that $\alpha'=-\alpha$.
  
The sufficiency of the condition in part (b) follows directly from
part (a) and Corollary 3.1. For the converse, it is enough to notice
that if for some $G\subset\{1,\ldots, n\}$, there is 
$0\neq u\in\H^{i-2}(\Delta_G; k)$, then there is $H\subset G$
such that $X^H u$ corresponds to a nonzero element in
$\bigcap_{j\in G\setminus H}\Ker(\H^{i-2}(\Delta_{G\setminus H}; k)
\longrightarrow\H^{i-2}(\Delta_{G\setminus(H\cup j)}; k))$. \Box  

\smallskip

\proclaim Example 1. Let $R=k[a,b,c,d]$ and 
$B=(ab, bc, cd, ad, ac)$. Then $\Delta$ is the simplicial complex:
$$\epsfbox{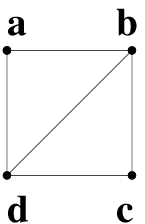}$$

Theorem 3.4(a) gives easily that
$${\rm Ass}(\Ext^3_R(R/B, R))=\{(a,b,d), (b,c,d)\}.$$

\proclaim Example 2. In general, it is not sufficient for 
$\beta_{|\alpha|-i, \alpha}(B^{\vee})$ to be nonzero in order to have
$P_{\alpha}\in {\rm Ass}(\Ext^i_R(R/B, R))$.

Let's consider $R=k[a,b,c]$ and $B=(a,bc)$. Then $\Delta$
is the simplicial complex:
$$\epsfbox{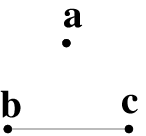}$$

Using Theorem 3.4(a), we get:
$${\rm Ass}(\Ext^2_R(R/B, R))=\{(a,b), (a,c)\},$$
while
$$\{F\,\vert\,\beta_{|\alpha_F|-2, \alpha_F}(B^{\vee})\neq 0\}=
\{ \{a,b,c\}, \{a,b\}, \{a,c\}\}.$$

\bigskip

\bigskip

\centerline{\bf References}
D.Bayer, H.Charalambous, and S.Popescu (1998).
 {\sl Extremal Betti numbers
  and applications to monomial ideals}, preprint.

W.Bruns and J.Herzog (1993). {\sl Cohen Macaulay rings}, Cambridge
   Univ. Press.

J.Eagon and V.Reiner (1996). {\sl Resolutions of Stanley-Reisner rings
   and Alexander duality}, preprint. 

D.Eisenbud (1995). {\sl Commutative Algebra with a View Toward
   Algebraic Geometry}, Springer.

D.Eisenbud, M.Musta\c{t}\v{a} and M.Stillman (1998). {\sl Cohomology
   on toric varieties and local cohomology with
monomial support}, preprint.

R.Godement (1958). {\sl Topologie algebrique et theorie des faisceaux},
   Paris, Herman.
  
M.Hochster (1977). {\sl Cohen-Macaulay rings, combinatorics and
   simplicial complexes}, in Ring theory II, B.R.McDonald, R.A.Morris (eds),
   Lecture Notes in Pure and Appl. Math., 26, M.Dekker.

J.R.Munkres (1984). {\sl Elements of algebraic topology},
    Benjamin/Cummings, Menlo Park CA.

M.Musta\c{t}\v{a} (1999). {\sl Vanishing theorems on toric varieties},
   in preparation.

N.Terai (1997). {\sl Generalization of Eagon-Reiner theorem and
   h-vectors of graded rings}, preprint.

N.Terai (1998). {\sl Local cohomology with respect to monomial
   ideals}, in preparation.

W.V.Vasconcelos (1998). {\sl Computational methods in commutative algebra
   and algebraic geometry}, Algorithms and Computation in Mathematics, vol.2,
   Springer--Verlag.

U.Walther (1999). {\sl Algorithmic computation of local cohomology modules
   and the cohomological dimension of algebraic varieties},
   Journal for Pure and Applied Algebra, to appear.

K.Yanagawa (1998). {\sl Alexander duality for 
   Stanley-Reisner rings and squarefree} 
   $N^n-$  {\sl graded modules},
   preprint.
 
\bigskip
\vbox
{\noindent Author Adress:} 
\smallskip
{\noindent Mircea Mustata}\par
{\noindent Department of Mathematics, Univ. of California, Berkeley;
Berkeley CA 94720}\par
{\noindent mustata@math.berkeley.edu}
\smallskip
{\noindent Institute of Mathematics of the Romanian Academy,
Calea Grivitei 21, Bucharest, Romania}\par
{\noindent mustata@stoilow.imar.ro}

\end